# SELECTING OPTIMAL MULTISTEP PREDICTORS FOR AUTOREGRESSIVE PROCESSES OF UNKNOWN ORDER

By Ching-Kang Ing

*Academia Sinica and National Taiwan University*

We consider the problem of choosing the optimal (in the sense of mean-squared prediction error) multistep predictor for an autoregressive (AR) process of finite but unknown order. If a working AR model (which is possibly misspecified) is adopted for multistep predictions, then two competing types of multistep predictors (i.e., plug-in and direct predictors) can be obtained from this model. We provide some interesting examples to show that when both plug-in and direct predictors are considered, the optimal multistep prediction results cannot be guaranteed by correctly identifying the underlying model's order. This finding challenges the traditional model (order) selection criteria, which usually aim to choose the order of the true model. A new prediction selection criterion, which attempts to seek the best combination of the prediction order and the prediction method, is proposed to rectify this difficulty. When the underlying model is stationary, the validity of the proposed criterion is justified theoretically. To obtain this result, asymptotic properties of accumulated squares of multistep prediction errors are investigated. In addition to overcoming the above difficulty, some other advantages of the proposed criterion are also mentioned.

**1. Introduction and overview.** In recent years there has been growing interest in the study of multistep prediction in various time series models [e.g., Findley (1984), Tiao and Xu (1993), Bhansali (1996, 1997), Haywood and Tunnicliffe-Wilson (1997), Hurvich and Tsai (1997), Findley, Pötscher and Wei (2001, 2003) and Ing (2003), among others]. Through these previous efforts, some new parameter estimation, prediction and model selection theories related to this research topic have been established. However, the problem of how to choose models to minimize multistep mean-squared pre-diction error (MSPE) has still not been clarified even for autoregressive (AR) processes. This motivated our study.

To fix ideas, let us assume that observations $x_1, \dots, x_n$ are generated from the stationary AR model

$$x_{t+1} = \sum_{i=1}^{p_1} a_i x_{t+1-i} + \varepsilon_{t+1}, \tag{1.1}$$

where $1 \le p_1 < \infty$ is unknown, $a_{p_1} \ne 0$, the $\varepsilon_t$'s are (unobservable) uncorrelated random noises with zero mean and common variance $\sigma^2$, and the characteristic polynomial $A(z) = 1 - a_1 z - \cdots - a_{p_1} z^{p_1}$ has no zeros inside or on the unit circle. This last assumption implies that $x_{t+1}$ has a one-sided infinite moving-average representation

$$x_{t+1} = \sum_{i=0}^{\infty} b_i \varepsilon_{t+1-i},$$

where $b_i = 1$ for $i = 0$ and $|b_i| \le c_0 e^{-c_1 i}$ for $i \ge 1$ and some positive numbers $c_0$ and $c_1$. For later reference we also define the parameter space of interest:

$$\Lambda = \{(d_1, \dots, d_{p_1})' : -\infty < d_i < \infty \text{ for } 1 \le i \le p_1 \text{ and } \\ 1 - d_1 z - \cdots - d_{p_1} z^{p_1} \ne 0 \text{ for any complex number } |z| \le 1\}.$$

To predict $x_{n+h}$, $h \ge 1$, under the situation where $p_1$ is unknown, it is common to use a working AR model, which is possibly misspecified, to replace the true underlying AR($p_1$) model. Then a natural predictor of $x_{n+h}$ can be obtained by repeatedly using the fitted (by least squares) working model with the unknown future values replaced by their own forecasts. In the following discussion this predictor is referred to as the plug-in predictor. More specifically, let the order of the working AR model be denoted by $k$ and let the least-squares estimator of the coefficient vector in the working model be denoted by $\hat{\mathbf{a}}_n(1,k) = (\hat{a}_{1,n}(k), \dots, \hat{a}_{k,n}(k))'$, where $\hat{\mathbf{a}}_n(1,k)$ satisfies

$$\hat{\Gamma}_n(1,k)\hat{\mathbf{a}}_n(1,k) = \frac{1}{n-k} \sum_{j=k}^{n-1} \mathbf{x}_j(k) x_{j+1}$$

with $\mathbf{x}(k) = (x_j \dots, x_{j-k+1})'$ and

$$\hat{\Gamma}_n(h,k) = \frac{1}{n-h-k+1} \sum_{j=k}^{n-h} \mathbf{x}_j(k)\mathbf{x}_j'(k).$$

Then, for $h \ge 1$ the plug-in predictor can be expressed by

$$\hat{x}_{n+h}(k) = \mathbf{x}_n'(k)\hat{\mathbf{a}}_n(h,k), \tag{1.2}$$

where $\hat{\mathbf{a}}_n(h,k) = \hat{A}_n^{h-1}(k)\hat{\mathbf{a}}_n(1,k)$, and with $I_m$ and $\mathbf{0}_m$, respectively, denoting an identity matrix and a vector of zeros of dimension $m$,

$$\hat{A}_n(k) = \left(\hat{\mathbf{a}}_n(1,k) \;\middle|\; \frac{I_{k-1}}{\mathbf{0}'_{k-1}}\right).$$

(Note that $\hat{A}_n^0(k) = I_k$.) On the other hand, the direct predictor of $x_{n+h}$, $\check{x}_{n+h}(k)$, suggested by Findley (1984), is also frequently used as an alternative, where $\check{x}_{n+h}(k)$ is obtained through a linear least-squares regression of $x_{t+h}$ on $x_t, \dots, x_{t-k+1}$, that is,

$$\check{x}_{n+h}(k) = \mathbf{x}'_n(k)\check{\mathbf{a}}_n(h,k), \tag{1.3}$$

where $\check{\mathbf{a}}_n(h,k)$ satisfies

$$\hat{\Gamma}_n(h,k)\check{\mathbf{a}}_n(h,k) = \frac{1}{n-h-k+1}\sum_{j=k}^{n-h}\mathbf{x}_j(k)x_{j+h}.$$

Viewing (1.2) and (1.3), it is obvious that the plug-in and direct predictors are identical when $h = 1$. For $h \geq 2$ Ing [(2003), Theorems 1 and 2] showed that the plug-in predictor has an advantage over the direct predictor in situations where the order of the working model, $k$, is not less than $p_1$. More specifically, as $h \geq 2$ and $k \geq p_1$, the MSPE of the plug-in predictor,

$$\text{MSPE}\, P_{n,h}(k) = E(x_{n+h} - \hat{x}_{n+h}(k))^2,$$

and that of the direct predictor,

$$\text{MSPE}\, D_{n,h}(k) = E(x_{n+h} - \check{x}_{n+h}(k))^2,$$

have the property

$$\lim_{n\to\infty} \frac{\text{MSPE}\, D_{n,h}(k) - \sigma_h^2}{\text{MSPE}\, P_{n,h}(k) - \sigma_h^2} > 1, \tag{1.4}$$

where $\sigma_h^2 = \sigma^2 \sum_{j=0}^{h-1} b_j^2$. Therefore, $\hat{x}_{n+h}(k)$ is asymptotically more efficient than $\check{x}_{n+h}(k)$ when $k \geq p_1$ and $h \geq 2$. For more details, see (2.2)–(2.4) of Section 2. Ing (2003) also compared the prediction efficiencies of $\hat{x}_{n+h}(k)$ and $\hat{x}_{n+h}(k+1)$ and those of $\check{x}_{n+h}(k)$ and $\check{x}_{n+h}(k+1)$ for $k \geq p_1$. Under certain conditions it was shown in Theorem 3 of Ing (2003) (see also Theorem 2.3 of Section 2) that

$$\lim_{n\to\infty} \frac{\text{MSPE}\, P_{n,h}(k+1) - \sigma_h^2}{\text{MSPE}\, P_{n,h}(k) - \sigma_h^2} > 1 \tag{1.5}$$

and

$$\lim_{n\to\infty} \frac{\text{MSPE}\, D_{n,h}(k+1) - \sigma_h^2}{\text{MSPE}\, D_{n,h}(k) - \sigma_h^2} > 1 \tag{1.6}$$

hold for $h \geq 1$ and $k \geq p_1$. Inequalities (1.4)–(1.6) suggest that from the MSPE point of view, $\hat{x}_{n+h}(p_1)$ seems to be the optimal choice among two competing families of candidate predictors,

$$\text{family I} = \{\hat{x}_{n+h}(1), \ldots, \hat{x}_{n+h}(K)\}$$

and

$$\text{family II} = \{\check{x}_{n+h}(1), \ldots, \check{x}_{n+h}(K)\},$$

where $K$ is known to satisfy $K \geq p_1$. [Note that we sometimes use $(k,1)$ to denote $\hat{x}_{n+h}(k)$ and use $(k,2)$ to denote $\check{x}_{n+h}(k)$.] Surprisingly, when $h \geq 2$ this conjecture is not true, provided $(a_1, \ldots, a_{p_1})'$ falls into some nonempty subset of $\Lambda$.

To see this, let us begin with the linear predictor of $x_{t+h}, h \geq 1$, based on the infinite past, $x_{t-j}, j \geq 0$, with the smallest MSPE. Let this predictor be denoted by $\tilde{x}_{t+h}$. Then we have

$$\tilde{x}_{t+h} = \sum_{j=1}^{p_h} a_j(h, p_h) x_{t+1-j},$$

where $a_{p_h}(h, p_h) \neq 0$ and

$$(a_1(h, p_h), \ldots, a_{p_h}(h, p_h))' = \mathbf{a}_D(h, p_h)$$

with $\mathbf{a}_D(h,k) = \Gamma^{-1}(k)(\gamma_h, \ldots, \gamma_{h+k-1})'$, $\Gamma(k) = E(\mathbf{x}_1(k)\mathbf{x}_1'(k))$ and $\gamma_j = E(x_t x_{t-j})$. We also have

$$x_{t+h} = \tilde{x}_{t+h} + \eta_{t,h}, \tag{1.7}$$

where $\eta_{t,h} = \sum_{j=0}^{h-1} b_j \varepsilon_{t+h-j}$. Model (1.7) is referred to as the $h$-step prediction model that corresponds to model (1.1) [note that when $h = 1$, $a_j(1, p_1) = a_j$ for $j = 1, \ldots, p_1$]. One notable but often disregarded feature of model (1.7) is that when $h > 1$, $p_h$ can be strictly less than $p_1$ and vary with $h$. For example, if $p_1 = 2$, then the corresponding two-step prediction model is

$$x_{t+2} = (a_1^2 + a_2)x_t + a_2 a_1 x_{t-1} + \varepsilon_{t+2} + a_1 \varepsilon_{t+1}.$$

Hence $p_2 = 1 < p_1$ if $a_1 = 0$. A similar situation also arises in the three-step prediction case, provided that $a_1^2 + a_2 = 0$. This phenomenon can occur even if all parameters in the one-step prediction model are large in magnitude. This also creates some unexpected difficulties in assessing the performances of the plug-in and direct predictors.

Note that when $p_h < p_1$ it seems more interesting to compare the performances of $\hat{x}_{n+h}(p_1)$ and $\check{x}_{n+h}(p_h)$ rather than those of $\hat{x}_{n+h}(k)$ and $\check{x}_{n+h}(k)$. In Section 2, some interesting examples are given to show that when $p_h < p_1$ and $h \geq 2$,

$$\lim_{n \to \infty} \frac{\text{MSPE}\, D_{n,h}(p_h) - \sigma_h^2}{\text{MSPE}\, P_{n,h}(p_1) - \sigma_h^2} < 1 \tag{1.8}$$

can occur. Moreover, since the value of the above limit depends on unknown parameters, it is not possible to determine the rankings of $\hat{x}_{n+h}(p_1)$ and $\check{x}_{n+h}(p_h)$ from the point of view of MSPE. This phenomenon further leads us to face a fundamental problem while selecting multistep predictors; that is, instead of the multistep predictor obtained by identifying the one-step prediction model's order, can a multistep predictor be constructed to minimize the multistep MSPE directly? As mentioned, this problem is complicated when both families I and II are considered. In this situation, the prediction order and the prediction method must be taken into account simultaneously.

This article aims to resolve the above problem. The strategy adopted herein is to find a statistic for each $\mathrm{MSPE}\,P_{n,h}(k)$ and $\mathrm{MSPE}\,P_{n,h}(k)$, $k = 1, \ldots, K$, and to show that the ordering of these statistics coincides with the ordering of their corresponding multistep MSPEs. To achieve this goal, we consider the multistep generalizations of accumulated prediction errors (APEs) based on sequential plug-in and direct predictors, namely,

$$\mathrm{APE}\,P_{n,h}(k) = \sum_{i=m_h}^{n-h} (x_{i+h} - \hat{x}_{i+h}(k))^2 \tag{1.9}$$

and

$$\mathrm{APE}\,D_{n,h}(k) = \sum_{i=m_h}^{n-h} (x_{i+h} - \check{x}_{i+h}(k))^2, \tag{1.10}$$

respectively, where $m_h$ denotes the smallest positive number such that $\hat{\mathbf{a}}_i(h, K)$ and $\check{\mathbf{a}}_i(h, K)$ are well defined for all $i \geq m_h$. Note that the APE with $h = 1$, namely, $\mathrm{APE}\,P_{n,1}(k) = \mathrm{APE}\,D_{n,1}(k)$, was first proposed by Rissanen (1986) for the purpose of determining $p_1$. Subsequently, the statistical properties of $\mathrm{APE}\,P_{n,1}(k)$ were investigated by Wei (1987, 1992) in stochastic regression models, which included model (1.1) as a special case. However, as indicated in Section 3, Wei's approach cannot be directly applied to the case of $h \geq 2$. Theorems 3.1 and 3.2 (also in Section 3) are devoted to dealing with this difficulty. In particular, the results obtained in these theorems show that the ordering of the multistep MSPEs of the predictors in families I and II can be well preserved by their corresponding multistep APEs when $n$ is sufficiently large. Based on this finding, we propose the following predictor selection procedure $(\hat{k}_n, \hat{j}_n)$, where $1 \leq \hat{k}_n \leq K$ and $1 \leq \hat{j}_n \leq 2$ (recall that $\hat{k}_n$ denotes the prediction order and $\hat{j}_n$ denotes the method of prediction):

STEP 1. Define $\hat{k}^{(1)}_{D,n} = \arg\min_{1 \leq k \leq K} \mathrm{APE}\,D_{n,1}(k)$.

STEP 2. Define

$$\hat{k}_{D,n}^{(h)} = \arg \min_{1 \le k \le K} \mathrm{APE}\, D_{n,h}(k)$$

and define

$$\hat{k}_n^{(1,h)} = \arg \min_{\hat{k}_{D,n}^{(1)} \le k \le K} \mathrm{APE}\, P_{n,h}(k).$$

STEP 3. If $\mathrm{APE}\, D_{n,h}(\hat{k}_{D,n}^{(h)}) > \mathrm{APE}\, P_{n,h}(\hat{k}_n^{(1,h)})$, then $(\hat{k}_n, \hat{j}_n) = (\hat{k}_n^{(1,h)}, 1)$; otherwise $(\hat{k}_n, \hat{j}_n) = (\hat{k}_{D,n}^{(h)}, 2)$.

We show in Theorem 3.4 of Section 3 that with probability 1, $(\hat{k}_n, \hat{j}_n)$ ultimately can choose the best predictor among families I and II regardless of whether $p_h < p_1$ or $p_h = p_1$. This property is referred to as the asymptotic efficiency; see Section 3 for the explicit definition. Moreover, $p_1$ can also be consistently estimated by $\hat{k}_{D,n}^{(1)}$.

It is worth noting that in this article more than a treatment of the difficulty caused by (1.8) is offered: (1) To the author's knowledge, $(\hat{k}_n, \hat{j}_n)$ is the first criterion that is designed to choose the optimal multistep predictor from the "honest" prediction point of view. By honest prediction, we mean the prediction for the future of the observed time series; see Rissanen (1987, 1989) for details. In the context of time series, most model selection criteria for prediction are obtained or justified under the assumption that the processes used for estimation and for prediction are independent; see, for example, finite prediction error [FPE; Akaike (1969)], Akaike information criterion [AIC; Akaike (1974)] and $S_n(k)$ [Shibata (1980)]. However, this type of prediction, which differs from Rissanen's idea, does not seem to be natural for time series analysis; see also Ing and Wei (2004). Recently, Ing and Wei (2004) obtained optimality for honest predictions of AIC (referred to as same-realization predictions in their article) in stationary AR($\infty$) processes. However, because their main concern was the case of one-step predictions, they did not deal with the problem of choosing the optimal combination of prediction order and prediction method. (2) This article shows that accumulated squares of sequential prediction errors can be used to choose a good predictor even in certain nonstandard situations. The sequential prediction error of $\mathrm{APE}\, P_{n,h}(k)$ with $h \ge 2$ involves a nonlinear transformation of the one-step least-squares estimators. While the sequential prediction error of $\mathrm{APE}\, D_{n,h}(k)$ with $h \ge 2$ is directly obtained from ($h$-step) least squares, its martingale structure no longer exists [see the discussion after (3.6)]. These nonstandard situations, which are not encountered with the one-step APE, challenge the validity of the multistep generalizations of APE for model (predictor) selection. By establishing the asymptotic efficiency of $(\hat{k}_n, \hat{j}_n)$, we clarify this ambiguity.

This article is organized as follows. In Section 2, some preliminary results from Ing (2003) and some examples that motivated this work are introduced. The asymptotic efficiency of $(\hat{k}_n, \hat{j}_n)$ is established in Section 3. In Section 4, an extension of the proposed criterion to subset autoregressions is given. Concluding remarks are given in Section 5. Some technical results, which are useful for obtaining the APE $P_{n,h}(k)$ asymptotic expression with $k \geq p_1$ are provided in the Appendix.

**2. Preliminary results and motivating examples.** Throughout this section, it is assumed that in model (1.1) the $\varepsilon_t$'s are i.i.d. random variables with mean 0 and variance $\sigma^2 > 0$. We also assume that the distribution function of $\varepsilon_1$, $F(\cdot)$, has the property, for some positive numbers $\alpha$, $\eta$ and $M$,

$$|F(x) - F(y)| \leq M|x - y|^{\alpha}, \tag{2.1}$$

provided $|x - y| < \eta$. Theorems 2.1 and 2.2 provide asymptotic expressions for MSPE $P_{n,h}(k)$ and MSPE $D_{n,h}(k)$ with $k \geq p_1$, respectively. Their proofs can be found in Theorems 1 and 2 of Ing (2003).

THEOREM 2.1. *Assume that $\{x_t\}$ satisfies model* (1.1). *Also assume that $\{\varepsilon_t\}$ satisfies* (2.1) *and*

$$E(|\varepsilon_1|^{\theta_h}) < \infty,$$

*where $\theta_h = \max\{8, 2(h+1)\} + \delta$ for some $\delta > 0$. Then, for $k \geq p_1$ and $h \geq 1$,*

$$n(\mathrm{MSPE}\, P_{n,h}(k) - \sigma_h^2) = f_{1,h}(k) + O(n^{-1/2}), \tag{2.2}$$

*where $f_{1,h}(k) = \mathrm{tr}(\Gamma(k) L_h(k) \Gamma^{-1}(k) L'_h(k))\sigma^2$ with $L_h(k) = \sum_{j=0}^{h-1} b_j A^{h-1-j}(k)$,*

$$A(k) = \left( \mathbf{a}_D(1,k) \;\middle|\; \frac{I_{k-1}}{\mathbf{0}'_{k-1}} \right)$$

*and $A^0(k) = I_k$.*

THEOREM 2.2. *Let the assumptions of Theorem 2.1 hold, with $\theta_h$ replaced by $8 + \delta$ for some $\delta > 0$. Then, for $k \geq p_h$ and $h \geq 1$,*

$$n(\mathrm{MSPE}\, D_{n,h}(k) - \sigma_h^2) = f_{2,h}(k) + O(n^{-1/2}), \tag{2.3}$$

*where $f_{2,h}(k) = \mathrm{tr}\{\Gamma^{-1}(k)\,\mathrm{cov}(\sum_{j=0}^{h-1} b_j \mathbf{x}_j(k))\}\sigma^2$ and, for a random vector $\mathbf{y}$, $\mathrm{cov}(\mathbf{y}) = E\{(\mathbf{y} - E(\mathbf{y}))(\mathbf{y} - E(\mathbf{y}))'\}$.*

Bhansali [(1997), Proposition 3.2] showed that for $k \geq p_1 \geq 1$ and $h \geq 2$,

$$\frac{f_{2,h}(k)}{f_{1,h}(k)} > 1. \tag{2.4}$$

Therefore, if $k \geq p_1 \geq 1$ and $h \geq 2$, then $\hat{x}_{n+h}(k)$ is asymptotically more efficient than $\check{x}_{n+h}(k)$ in the sense of (1.4). For example, assume $h = 2$ and $k \geq p_1 \geq 1$. Then

$$f_{2,2}(k) = \{k + (k+2)a_1^2\}\sigma^2, \tag{2.5}$$

and

$$f_{1,2}(k) = \{(k+2)a_1^2 + k - 1 + a_k^2\}\sigma^2. \tag{2.6}$$

(Note that $|a_{p_1}| < 1$ and $a_k = 0$ for $k \geq p_1$.) Hence, for $k \geq p_1$,

$$\lim_{n\to\infty} \frac{\operatorname{MSPE} D_{n,2}(k) - \sigma_h^2}{\operatorname{MSPE} P_{n,2}(k) - \sigma_h^2} - 1 = \frac{1 - a_k^2}{(k+2)a_1^2 + k - 1 + a_k^2} > 0.$$

The following theorem shows that $f_{1,h}(k)$ and $f_{2,h}(k)$ with $k \geq p_1$ are strictly increasing functions of $k$.

THEOREM 2.3. (i) *Assume $h \geq 1$ and $k \geq p_1$. Then*

$$\frac{f_{1,h}(k+1)}{f_{1,h}(k)} > 1, \tag{2.7}$$

*provided*

$$b_{h-1} \neq 0 \tag{2.8}$$

*or*

$$\mathbf{l}^* \neq \mathbf{0}_{k+1}, \tag{2.9}$$

*where with the convention that $b_j = 0$ for $j < 0$, $\mathbf{l}^* = (\sum_{i=0}^{h-1} b_{h-1-k-i}b_i, \ldots, \sum_{i=0}^{h-1} b_{h-1-i}b_i)'$ is a $(k+1)$-dimensional vector.*

(ii) *Assume $h \geq 1$ and $k \geq p_1$. Then*

$$\frac{f_{2,h}(k+1)}{f_{2,h}(k)} > 1. \tag{2.10}$$

REMARK 1. A proof of Theorem 2.3 can be found in Ing [(2003), Theorem 3]. When $1 \leq h \leq 5$, it can be shown that either (2.8) or (2.9) holds for all $k \geq p_1$, and hence (2.7) holds without extra constraints on the parameter space. However, for general $h$ (especially when $h \gg k$), we are not able to establish (2.7) without conditions (2.8) or (2.9). For more details on these conditions, see Ing [(2003), Remark 2].

As immediate consequences of Theorems 2.1–2.3, we obtain (1.5) and (1.6). Inequalities (1.4)–(1.6) seem to suggest that

$$\lim_{n\to\infty}\frac{E(x_{n+h}-\hat{x}_{n+h}(p_1))^2-\sigma_h^2}{E(x_{n+h}-\tilde{x}_{n+h}(k))^2-\sigma_h^2}\le 1, \tag{2.11}$$

where $\tilde{x}_{n+h}(k)$ is any predictor in family I or II. However, as indicated by Remark 1, when $h$ is large, (2.7) cannot be guaranteed without (2.8) or (2.9). Therefore, it is not clear whether (2.11) still holds in the situation where both (2.8) and (2.9) are violated. Moreover, we will show that (2.11) can fail when $p_h < p_1$. To see this, let us begin with a simple extension of Theorem 2.2, which provides an asymptotic expression for $\mathrm{MSPE}\, D_{n,h}(k)$ with $k \ge p_h$.

Corollary 2.4. *Let the assumptions of Theorem* 2.2 *hold. Then* (2.3) *holds with* $k \ge p_h$ *and* $h \ge 1$.

Since Corollary 2.4 can be shown by an argument similar to that used to show Theorem 2.2, we omit the details. When $p_h < p_1$, it would be more interesting to compare

$$\lim_{n\to\infty}(\mathrm{MSPE}\, P_{n,h}(p_1)-\sigma_h^2) \quad \text{and} \quad \lim_{n\to\infty}(\mathrm{MSPE}\, D_{n,h}(p_h)-\sigma_h^2)$$

rather than

$$\lim_{n\to\infty}(\mathrm{MSPE}\, P_{n,h}(k)-\sigma_h^2) \quad \text{and} \quad \lim_{n\to\infty}(\mathrm{MSPE}\, D_{n,h}(k)-\sigma_h^2).$$

The following two examples show that the advantage of the plug-in predictor can vanish in this kind of comparison.

Example 1. Let $h = 2$ and $p_2 < p_1$. Then we see that $b_1 = a_1 = 0$ and $p_2 = p_1 - 1$. This fact and Corollary 2.4 yield that $f_{2,2}(p_1) - f_{2,2}(p_2) = \sigma^2$. On the other hand, by (2.5) and (2.6) we have $f_{2,2}(p_1) - f_{1,2}(p_1) = (1 - a_{p_1}^2)\sigma^2$. Therefore, $f_{1,2}(p_1) - f_{2,2}(p_2) = a_{p_1}^2\sigma^2 > 0$. As a result, we have, for $p_1 - p_2 = 1$,

$$\lim_{n\to\infty}\frac{\mathrm{MSPE}\, D_{n,2}(p_2)-\sigma_2^2}{\mathrm{MSPE}\, P_{n,2}(p_1)-\sigma_2^2}=\frac{f_{2,2}(p_2)}{f_{1,2}(p_1)}<1$$

and hence $\check{x}_{n+2}(p_2)$ is asymptotically more efficient than $\hat{x}_{n+2}(p_1)$ in this case.

For general $h$, the ratio of $f_{2,h}(p_h)/f_{1,h}(p_1)$ can be larger or smaller than 1, as shown in the following example.

EXAMPLE 2. First assume that $p_1 = 2$ and $h = 3$. By (2.5), (2.6) and the fact that when $k \geq p_1$,

$$f_{2,h+1}(k) - f_{1,h+1}(k) = f_{2,h}(k) - f_{1,h}(k) + \mathbf{e}_k' L_h(k)\Gamma^{-1}(k)L_h'(k)\mathbf{e}_k\sigma^4$$

[see Section 2 of Ing (2003)], where $\mathbf{e}_k' = (1, 0, \ldots, 0)$ is a $k$-dimensional vector,

$$f_{2,3}(2) - f_{1,3}(2) = (1 - a_2^2)\sigma^2 + \mathbf{e}_2' L_2(2)\Gamma^{-1}(2)L_2'(2)\mathbf{e}_2\sigma^4.$$

Some algebraic manipulations yield $\mathbf{e}_2' L_2(2)\Gamma^{-1}(2)L_2'(2)\mathbf{e}_2\sigma^4 = (1 + a_2)(1 - a_2 - 4a_1^2 a_2)\sigma^2$. Therefore

$$f_{2,3}(2) - f_{1,3}(2) = 2(1 + a_2)(1 - a_2 - 2a_1^2 a_2)\sigma^2. \tag{2.12}$$

Note that Bhansali [(1997), page 442] indicated that $f_{2,3}(2) - f_{1,3}(2) = (1 + a_2)(1 - a_2 - 2a_1^2 a_2)\sigma^2$. However, one can see that the leading constant 2 on the right-hand side of (2.12) is needed by examining a simple example which assumes that $-1 < a_1 < 1$ and $a_2 = 0$.

Now, assume $b_2 = a_1^2 + a_2 = 0$. Then $p_3 = 1 < 2 = p_1$ and, in view of (2.12),

$$f_{2,3}(2) - f_{1,3}(2) = 2(1 + a_2)(1 - a_2 + 2a_2^2)\sigma^2. \tag{2.13}$$

By Corollary 2.4,

$$f_{2,3}(1) = \frac{1 - 4a_2 + a_2^2}{1 - a_2}\sigma^2, \tag{2.14}$$

and

$$f_{2,3}(2) - f_{2,3}(1) = \left(1 - a_2 + \frac{2a_2^2}{1 - a_2}\right)\sigma^2. \tag{2.15}$$

According to (2.13)–(2.15),

$$\frac{f_{2,3}(p_3)}{f_{1,3}(p_1)} = \frac{f_{2,3}(1)}{f_{1,3}(2)} = \frac{1 - 4a_2 + a_2^2}{-4a_2 + 2a_2^2 - 2a_2^3 + 4a_2^4}. \tag{2.16}$$

Let the rational function on the right-hand side of (2.16) be denoted by $g(a_2)$ and let the unique solution of the equation $g(a_2) = 1$ with $-1 < a_2 < 0$ be denoted by $T$. Then it can be shown that $T \approx -0.54977$, $g(a_2) < 1$ if $-1 < a_2 < T$ and $g(a_2) > 1$ if $T < a_2 < 0$. Therefore, when $h \geq 3$ and $p_h < p_1$, it is not possible to determine the rankings of $\hat{x}_{n+h}(p_1)$ and $\check{x}_{n+h}(p_h)$ without knowledge of the AR parameters.

To illustrate the results obtained in Example 2, four AR(2) models,

$$x_t = 0.9x_{t-1} - 0.81x_{t-2} + \varepsilon_t, \tag{2.17}$$

$$x_t = 0.8x_{t-1} - 0.64x_{t-2} + \varepsilon_t, \tag{2.18}$$

$$x_t = 0.6x_{t-1} - 0.36x_{t-2} + \varepsilon_t \tag{2.19}$$

and

$$x_t = 0.5x_{t-1} - 0.25x_{t-2} + \varepsilon_t, \tag{2.20}$$

are considered in our simulation study, where $\varepsilon_t$'s are independent and identically $\mathcal{N}(0,1)$ distributed. The empirical estimates of $(\text{MSPE}\, D_{n,3}(1) - \sigma_3^2)/(\text{MSPE}\, P_{n,3}(2) - \sigma_3^2)$ for the above four models are obtained based on 20,000 replications for $n = 150$, 300, 500 and 1000. These empirical estimates and corresponding limiting values [given by (2.16)] are summarized in Table 1. One can see from these empirical results that $\check{x}_{n+3}(1)$ is more efficient than $\hat{x}_{n+3}(2)$ for models (2.17) and (2.18), and is less efficient than $\hat{x}_{n+3}(2)$ for the other two models. This conclusion coincides with that obtained from (2.16). In addition, the empirical estimates of $(\text{MSPE}\, D_{n,3}(1) - \sigma_3^2)/(\text{MSPE}\, P_{n,3}(2) - \sigma_3^2)$ are rather close to their corresponding limiting values even for $n = 150$.

As a conclusion, we note that when both the plug-in and direct predictors are taken into account, the optimal multistep prediction results cannot be guaranteed by correctly identifying $p_1$ or $p_h$. Hence, a predictor selection criterion that directly aims at the minimal MSPE (among those of the predictors in families I and II) is called for.

## 3. Main results.

Since we attempt to choose a candidate predictor among families I and II that has having the minimal MSPE (at least for large $n$), the loss functions of the candidate plug-in and direct predictors are defined as

$$L_{1,h}(k) = \begin{cases} \lim_{n\to\infty} n(\text{MSPE}\, P_{n,h}(k) - \sigma_h^2), & \text{if } p_1 \le k \le K, \\ \infty, & \text{if } k < p_1 \end{cases} \tag{3.1}$$

and

$$L_{2,h}(k) = \begin{cases} \lim_{n\to\infty} n(\text{MSPE}\, D_{n,h}(k) - \sigma_h^2), & \text{if } p_h \le k \le K, \\ \infty, & \text{if } k < p_h, \end{cases} \tag{3.2}$$

TABLE 1
*Simulation results for* $(\text{MSPE}\, D_{n,3}(1) - \sigma_3^2)/(\text{MSPE}\, P_{n,3}(2) - \sigma_3^2)$

| $n$ | Model (2.17) | (2.18) | (2.19) | (2.20) |
|---|---|---|---|---|
| 150 | 0.700 | 0.891 | 1.398 | 1.719 |
| 300 | 0.688 | 0.843 | 1.365 | 1.782 |
| 500 | 0.649 | 0.879 | 1.365 | 1.762 |
| 1000 | 0.673 | 0.872 | 1.379 | 1.761 |
| $f_{2,3}(1)/f_{1,3}(2)$ | 0.667 | 0.868 | 1.382 | 1.76 |

respectively, where the existence of the above limits is ensured by Theorems 2.1 and 2.2. To ensure the prediction loss due to underspecification is much larger than the loss due to overspecification, the loss function values of $(k,1)$ with $k<p_1$ and of $(k,2)$ with $k<p_h$ are set to $\infty$. A predictor selection criterion, $(\tilde{k}_n,\tilde{j}_n)$ with $1\le\tilde{k}_n\le K$ and $1\le\tilde{j}_n\le 2$, is said to be asymptotically efficient if

$$P((\tilde{k}_n,\tilde{j}_n)\in C_{h,K}\text{ eventually})=1, \tag{3.3}$$

where

$$C_{h,K}=\Big\{(k,j):1\le k\le K, 1\le j\le 2\text{ and } L_{j,h}(k)=\min_{1\le k_0\le K,1\le j_0\le 2}L_{j_0,h}(k_0)\Big\}.$$

Therefore, with probability 1 $(\tilde{k}_n,\tilde{j}_n)$ can ultimately choose a predictor having the minimal loss function value.

REMARK 2. Note that $C_{h,K}$ can contain more than one element. To see this, assume that $h=3$, $p_1=2$, $a_1^2+a_2=0$, $a_2=T\approx-0.54977$ and $K\ge 2$. (Recall that $p_3=1<p_1$ in this case.) By Theorems 2.1 and 2.3, Corollary 2.4 and Remark 1, we have $f_{1,3}(k)<f_{1,3}(k+1)$, $f_{2,3}(k)<f_{2,3}(k+1)$ and $f_{1,3}(k)<f_{2,3}(k)$ for $k\ge 2$. Moreover, by Example 2, $f_{1,3}(2)=f_{2,3}(1)$. As a result there are two elements, namely (1, 2) and (2, 1), in $C_{3,K}$.

The goal of this section is to show that (3.3) is fulfilled by $(\hat{k}_n,\hat{j}_n)$. We assume in this section that $\{\varepsilon_t\}$ in model (1.1) is a martingale difference sequence with respect to an increasing sequence of $\sigma$-fields $\{\mathcal{F}_t\}$, that is, $\varepsilon_t$ is $\mathcal{F}_t$-measurable, and $E(\varepsilon_t|\mathcal{F}_{t-1})=0$ a.s. for all $t$. We also assume that for some $\alpha>2$,

$$E(\varepsilon_t^2|\mathcal{F}_{t-1})=\sigma^2\quad\text{and}\quad\sup_t E(|\varepsilon_t|^\alpha|\mathcal{F}_{t-1})<\infty\qquad\text{a.s.} \tag{3.4}$$

Note that for $k\ge p_1$,

$$\operatorname{APE}P_{n,h}(k)=\sum_{i=m_h}^{n-h}\{\eta_{i,h}-\mathbf{x}_i'(k)\hat{L}_{i,h}(k)(\hat{\mathbf{a}}_i(1,k)-\mathbf{a}_D(1,k))\}^2 \tag{3.5}$$

and for $k\ge p_h$,

$$\operatorname{APE}D_{n,h}(k)=\sum_{i=m_h}^{n-h}\{\eta_{i,h}-\mathbf{x}_i'(k)(\check{\mathbf{a}}_i(h,k)-\mathbf{a}_D(h,k))\}^2, \tag{3.6}$$

where $\eta_{i,h}$ is defined in (1.7) and $\hat{L}_{i,h}(k)=\sum_{j=0}^{h-1}b_j\hat{A}_i^{h-1-j}(k)$, with $\hat{A}_i^{h-1-j}(k)$ defined below (1.2). The asymptotic properties of $\operatorname{APE}P_{n,h}(k)=\operatorname{APE}D_{n,h}(k)$

with $h=1$ were investigated by Wei (1987, 1992) in stochastic regression models. One of the key steps in Wei's analysis is to express the (second-order) residual sum of squares of the fitted (by least squares) model in a recursive form. His approach, however, cannot be directly applied to the situation considered in this article. This is because for $\operatorname{APE} P_{n,h}(k)$ with $h \geq 2$ there is a random matrix $\hat{L}_{i,h}(k)$ that lies between $\mathbf{x}_i'(k)$ and $(\hat{\mathbf{a}}_i(1,k) - \mathbf{a}_D(1,k))$, and for $\operatorname{APE} D_{n,h}(k)$ with $h \geq 2$ the rightmost component $\sum_{j=k}^{n-h} \mathbf{x}_j(k)\eta_{j,h}$ of the centered estimator

$$\check{\mathbf{a}}_i(h,k) - \mathbf{a}_D(h,k) = \frac{1}{i-k-h+1}\hat{\Gamma}_i^{-1}(h,k)\sum_{j=k}^{n-h}\mathbf{x}_j(k)\eta_{j,h}$$

is no longer a martingale transformation. Therefore, some new technical tools are needed to overcome these difficulties.

Theorems 3.1 and 3.2 describe the asymptotic behavior of $\operatorname{APE} P_{n,h}(k)$ and $\operatorname{APE} D_{n,h}(k)$ in the correctly specified case.

THEOREM 3.1. *Assume that* $\{x_t\}$ *satisfies model* (1.1). *Also assume condition* (3.4). *Then for* $k \geq p_1$ *and* $h \geq 1$,

$$\operatorname{APE} P_{n,h}(k) - \sum_{i=m_h}^{n-h} \eta_{i,h}^2 = \sigma^2 f_{1,h}(k)\log n + o(\log n) \qquad a.s. \tag{3.7}$$

PROOF. Rewrite the right-hand side of (3.5) as

$$\sum_{i=m_h}^{n-h}(\eta_{i,h})^2 - 2\sum_{i=m_h}^{n-h}\{\mathbf{x}_i'(k)\hat{L}_{i,h}(k)(\hat{\mathbf{a}}_i(1,k)-\mathbf{a}_D(1,k))\}\eta_{i,h}$$
$$+ \sum_{i=m_h}^{n-h}\{\mathbf{x}_i'(k)\hat{L}_{i,h}(k)(\hat{\mathbf{a}}_i(1,k)-\mathbf{a}_D(1,k))\}^2.$$

This and Chow (1965) yield that

$$\begin{aligned}\operatorname{APE} P_{n,h}(k) &- \sum_{i=m_h}^{n-h}(\eta_{i,h})^2 \\ &= \sum_{i=m_h}^{n-h}\{\mathbf{x}_i'(k)\hat{L}_{i,h}(k)(\hat{\mathbf{a}}_i(1,k)-\mathbf{a}_D(1,k))\}^2(1+o(1)) + O(1) \qquad \text{a.s.}\end{aligned} \tag{3.8}$$

To deal with the right-hand side of (3.8), we first introduce $Q_n^*(h,k)$, where

$$Q_n^*(h,k) = \left(\sum_{j=k}^{n-h}\mathbf{x}_j(k)\varepsilon_{j+1}\right)' S' V_{n-h} S \left(\sum_{j=k}^{n-h}\mathbf{x}_j(k)\varepsilon_{j+1}\right) \tag{3.9}$$

with $S=\Gamma(k)L_h(k)\Gamma^{-1}(k)$ and $V_i=(\sum_{j=k}^{i}\mathbf{x}_j(k)\mathbf{x}_j'(k))^{-1}$.

Following Lai and Wei [(1982), equation (2.16)], we obtain a recursive expression for $Q_n^*(h,k)$,

$$\begin{aligned}
&Q_n^*(h,k)+\sum_{i=m_h}^{n-h}\left\{\mathbf{x}_i'(k)V_{i-1}S\left(\sum_{j=k}^{i-1}\mathbf{x}_j(k)\varepsilon_{j+1}\right)\right\}^2 c_i^{-1}\\
&\qquad =Q_{m_h+h-1}^*(h,k)+\sum_{i=m_h}^{n-h}\mathbf{x}_i'(k)S'V_{i-1}S\mathbf{x}_i(k)\varepsilon_{i+1}^2\\
&\qquad\quad +\mathrm{I}+\mathrm{II}+\mathrm{III},
\end{aligned}\tag{3.10}$$

where

$$\mathrm{I}=2\sum_{i=m_h}^{n-h}\mathbf{x}_i'(k)S'V_{i-1}S\left(\sum_{j=k}^{i-1}\mathbf{x}_j(k)\varepsilon_{j+1}\right)\varepsilon_{i+1},$$

$$\mathrm{II}=-2\sum_{i=m_h}^{n-h}(\mathbf{x}_i'(k)S'V_{i-1}\mathbf{x}_i(k))\mathbf{x}_i'(k)V_{i-1}S\left(\sum_{j=k}^{i-1}\mathbf{x}_j(k)\varepsilon_{j+1}\right)\varepsilon_{i+1}c_i^{-1}$$

and

$$\mathrm{III}=-\sum_{i=m_h}^{n-h}(\mathbf{x}_i'(k)S'V_{i-1}\mathbf{x}_i(k))^2\varepsilon_{i+1}^2c_i^{-1}$$

with $c_i=(1+\mathbf{x}_i'(k)V_{i-1}\mathbf{x}_i(k))$. By (3.4), Theorem 2 of Lai and Wei (1985) and the martingale strong law of Lai and Wei (1982), we have

$$\lim_{n\to\infty}\frac{1}{n}V_n^{-1}=\Gamma(k)\qquad\text{a.s.},\tag{3.11}$$

which together with (3.4) and an analogy with (2.31) of Wei (1987) yields

$$Q_n^*(h,k)=o(\log n)\qquad\text{a.s.}\tag{3.12}$$

Since $c_n=(1-\mathbf{x}_n'(k)V_n\mathbf{x}_n(k))^{-1}$, by Theorem 4 of Lai and Wei (1983) [which ensures that $\lim_{n\to\infty}\mathbf{x}_n'(k)V_n\mathbf{x}_n(k)=0$ a.s.], we have

$$\lim_{n\to\infty}c_n=1\qquad\text{a.s.}\tag{3.13}$$

Now, by (3.4), (3.12), (3.13) and Chow (1965), we can rewrite (3.10) as

$$\begin{aligned}
&(1+o(1))\sum_{i=m_h}^{n-h}\left\{\mathbf{x}_i'(k)V_{i-1}S\left(\sum_{j=k}^{i-1}\mathbf{x}_j(k)\varepsilon_{j+1}\right)\right\}^2\\
&\qquad =o(\log n)+O(1)+(1+o(1))\sigma^2\sum_{i=m_h}^{n-h}\mathbf{x}_i'(k)S'V_{i-1}S\mathbf{x}_i(k)\\
&\qquad\quad +\mathrm{I}+\mathrm{II}+\mathrm{III}\qquad\text{a.s.}
\end{aligned}\tag{3.14}$$

Reasoning as in the proof of Lemma 2.1 of Wei (1992), we obtain

$$\lim_{n\to\infty}\frac{\sigma^2}{\log n}\sum_{i=m_h}^{n-h}\mathbf{x}_i'(k)S'V_{i-1}S\mathbf{x}_i(k)=f_{1,h}(k) \qquad \text{a.s.} \tag{3.15}$$

It is shown in the Appendix that

$$\mathrm{I}=o(\log n) \qquad \text{a.s.} \quad \text{and} \quad \mathrm{II}=o(\log n) \qquad \text{a.s.} \tag{3.16}$$

Moreover, by (3.11), Theorem 3 of Lai and Wei (1983), (2.10) and (2.12) of Lai and Wei (1982), and an analogy with Lemma 2.1 of Wei (1992),

$$\begin{aligned}\mathrm{III}&=O(1)+o\left(\sum_{i=m_h}^{n-h}|\mathbf{x}_i'(k)S'V_{i-1}\mathbf{x}_i'(k)|\varepsilon_{i+1}^2\right) \qquad \text{a.s.}\\ &=o(\log n) \qquad \text{a.s.}\end{aligned}$$

This, together with (3.14)–(3.16), yields

$$\begin{aligned}&\sum_{i=m_h}^{n-h}\left\{\mathbf{x}_i'(k)V_{i-1}S\left(\sum_{j=k}^{i-1}\mathbf{x}_j(k)\varepsilon_{j+1}\right)\right\}^2\\ &\qquad=\sigma^2 f_{1,h}(k)\log n+o(\log n) \qquad \text{a.s.}\end{aligned} \tag{3.17}$$

In view of (3.8) and (3.17) this proof is completed if we can show that

$$\begin{aligned}&\sum_{i=m_h}^{n-h}\{\mathbf{x}_i'(k)\hat{L}_{i,h}(k)(\hat{\mathbf{a}}_i(1,k)-\mathbf{a}_D(1,k))\}^2\\ &\qquad=\sum_{i=m_h}^{n-h}\left\{\mathbf{x}_i'(k)V_{i-1}\hat{S}_i\left(\sum_{j=k}^{i-1}\mathbf{x}_j(k)\varepsilon_{j+1}\right)\right\}^2\\ &\qquad=\sum_{i=m_h}^{n-h}\left\{\mathbf{x}_i'(k)V_{i-1}S\left(\sum_{j=k}^{i-1}\mathbf{x}_j(k)\varepsilon_{j+1}\right)\right\}^2\\ &\qquad\quad+o(\log n) \qquad \text{a.s.,}\end{aligned} \tag{3.18}$$

where $\hat{S}_i=V_{i-1}^{-1}\hat{L}_{i,h}(k)V_{i-1}$. Since by (3.11) and Theorem 1 of Lai and Wei (1983) $\lim_{n\to\infty}\hat{S}_n=S$ a.s., this fact and (A.1) imply that

$$\sum_{i=m_h}^{n-h}\left\{\mathbf{x}_i'(k)V_{i-1}(\hat{S}_i-S)\left(\sum_{j=k}^{i-1}\mathbf{x}_j(k)\varepsilon_{j+1}\right)\right\}^2=o(\log n) \qquad \text{a.s.} \tag{3.19}$$

Consequently, (3.18) follows from (3.17), (3.19) and the Cauchy–Schwarz inequality. □

THEOREM 3.2. *Let the assumptions of Theorem* 3.1 *hold. Then for* $k \geq p_h$ *and* $h \geq 1$,

$$\text{APE}\, D_{n,h}(k) - \sum_{i=m_h}^{n-h} \eta_{i,h}^2 = \sigma^2 f_{2,h}(k) \log n + o(\log n) \qquad a.s. \tag{3.20}$$

PROOF. We only show (3.20) for $h = 2$, because the result for $h \geq 3$ can be obtained similarly and that for $h = 1$ was verified in Wei (1992). Reasoning as for (3.8), we have, for $k \geq p_h$,

$$\begin{aligned} &\text{APE}\, D_{n,2}(k) - \sum_{i=m_h}^{n-2} (\eta_{i,2})^2 \\ &\quad = (1 + o(1)) \sum_{i=m_h}^{n-2} \left\{ \mathbf{x}_i'(k) V_{i-2} \left( \sum_{j=k}^{i-2} \mathbf{x}_j(k) \eta_{j,2} \right) \right\}^2 + O(1) \qquad \text{a.s.} \end{aligned} \tag{3.21}$$

Now consider

$$\bar{Q}_n(2,k) = \left( \sum_{i=k}^{n-2} \mathbf{x}_i(k) \eta_{i,2} \right)' V_{n-2} \left( \sum_{i=k}^{n-2} \mathbf{x}_i(k) \eta_{i,2} \right).$$

Following Theorem 1 of Wei (1987) and (3.14), we have

$$(1 + o(1)) T(k) = \bar{Q}_{m_h+1}(2,k) - \bar{Q}_n(2,k) + B(k) + C(k), \tag{3.22}$$

where

$$T(k) = \sum_{i=m_h}^{n-2} \left\{ \mathbf{x}_i'(k) V_{i-1} \left( \sum_{j=k}^{i-1} \mathbf{x}_j(k) \eta_{j,2} \right) \right\}^2,$$

$$B(k) = \sum_{i=m_h}^{n-2} \mathbf{x}_i'(k) V_i \mathbf{x}_i(k) \eta_{i,2}^2$$

and

$$C(k) = 2 \sum_{i=m_h}^{n-2} \mathbf{x}_i'(k) V_{i-1} \left( \sum_{j=k}^{i-1} \mathbf{x}_j(k) \eta_{j,2} \right) c_i^{-1} \eta_{i,2}.$$

[Notice that by Theorem 3 of Lai and Wei (1983) and (3.11), (3.13) still holds with $p_h \leq k < p_1$.]

In what follows we deal with $\bar{Q}_n(2,k)$, $B(k)$ and $C(k)$ separately. For $\bar{Q}_n(2,k)$, by an analogy with Theorem 3 of Wei (1987),

$$\begin{aligned} \bar{Q}_n(2,k) &= o\left( \log\left( \sum_{i=k-1}^{n-2} \|\mathbf{x}_i(k)\|^2 + \|a_1 \mathbf{x}_{i+1}(k)\|^2 \right) \right) \qquad \text{a.s.} \\ &= o(\log n) \qquad \text{a.s.}, \end{aligned} \tag{3.23}$$

where the second equality is ensured by (3.11).

For $B(k)$ we have

$$\begin{aligned} B(k) &= \sum_{i=m_h}^{n-2} \mathbf{x}_i'(k) V_i \mathbf{x}_i(k) \varepsilon_{i+2}^2 + a_1^2 \sum_{i=m_h}^{n-2} \mathbf{x}_i'(k) V_i \mathbf{x}_i(k) \varepsilon_{i+1}^2 \\ &\quad + 2a_1 \sum_{i=m_h}^{n-2} \mathbf{x}_i'(k) V_i \mathbf{x}_i(k) \varepsilon_{i+1}\varepsilon_{i+2}. \end{aligned} \tag{3.24}$$

According to Theorem 1 of Wei (1987), (3.11), (3.13) and Chow (1965), the right-hand side of (3.24) can be further expressed as

$$\begin{aligned} &\sigma^2(1+a_1^2)k\log n + o\left(\sum_{i=m_h}^{n-2} (\mathbf{x}_i'(k) V_i \mathbf{x}_i(k))^2 \varepsilon_{i+2}^2\right) + o(\log n) \qquad \text{a.s.} \\ &\quad = \sigma^2(1+a_1^2)k\log n + o(\log n) \qquad \text{a.s.} \end{aligned} \tag{3.25}$$

Therefore

$$B(k) = \sigma^2(1+a_1^2)k\log n + o(\log n) \qquad \text{a.s.} \tag{3.26}$$

To deal with $C(k)$, we have

$$\tfrac{1}{2}C(k) = D(k) + E(k) + F(k) + G(k) + H(k), \tag{3.27}$$

where

$$D(k) = \sum_{i=m_h}^{n-2} \mathbf{x}_i'(k) V_{i-1} \left(\sum_{j=k}^{i-2} \mathbf{x}_j(k)\eta_{j,2}\right) c_i^{-1}(\varepsilon_{i+2} + a_1\varepsilon_{i+1}),$$

$$E(k) = a_1^2 \sum_{i=m_h}^{n-2} \mathbf{x}_i'(k) V_{i-1} \mathbf{x}_{i-1}(k) c_i^{-1} \varepsilon_i \varepsilon_{i+1},$$

$$F(k) = a_1 \sum_{i=m_h}^{n-2} \mathbf{x}_i'(k) V_{i-1} \mathbf{x}_{i-1}(k) c_i^{-1} \varepsilon_{i+1}^2,$$

$$G(k) = \sum_{i=m_h}^{n-2} \mathbf{x}_i'(k) V_{i-1} \mathbf{x}_{i-1}(k) c_i^{-1} \varepsilon_{i+1} \varepsilon_{i+2},$$

$$H(k) = a_1 \sum_{i=m_h}^{n-2} \mathbf{x}_i'(k) V_{i-1} \mathbf{x}_{i-1}(k) c_i^{-1} \varepsilon_i \varepsilon_{i+2}.$$

By (3.4), (3.13) and Lemma 2(iii) of Lai and Wei (1982), we can show that

$$D(k) = o\left(\sum_{i=m_h}^{n-2} \left\{\mathbf{x}_i'(k) V_{i-1} \left(\sum_{j=k}^{i-2} \mathbf{x}_j(k)\eta_{j,2}\right)\right\}^2\right) + O(1) \qquad \text{a.s.} \tag{3.28}$$

Similarly,

$$
\begin{aligned}
E(k) &= o\left(\sum_{i=m_h}^{n-2}(\mathbf{x}_i'(k)V_{i-1}\mathbf{x}_{i-1}(k))^2\varepsilon_i^2\right)+O(1) \qquad \text{a.s.}\\
&= o\left(\sum_{i=m_h}^{n-2}\mathbf{x}_{i-1}'(k)V_{i-1}\mathbf{x}_{i-1}(k)\varepsilon_i^2\right)+O(1) \qquad \text{a.s.}\\
&= o(\log n) \qquad \text{a.s.},
\end{aligned} \tag{3.29}
$$

where the second equality is ensured by (3.13) and the Cauchy–Schwarz inequality, and the last equality is guaranteed by the same argument used to obtain Theorem 1 of Wei (1987). The same reasoning that shows (3.29) also gives

$$
G(k) = o(\log n) \qquad \text{a.s.} \tag{3.30}
$$

and

$$
H(k) = o(\log n) \qquad \text{a.s.} \tag{3.31}
$$

We now deal with $F(k)$. By an analogy with Lai and Wei (1982) we can show that

$$
\begin{aligned}
&\sum_{i=m_h}^{n-2}\mathbf{x}_i'(k)V_{i-1}\mathbf{x}_{i-1}(k)c_i^{-1}\varepsilon_{i+1}^2\\
&\quad= \sigma^2\sum_{i=m_h}^{n-2}\mathbf{x}_i'(k)V_{i-1}\mathbf{x}_{i-1}(k)c_i^{-1}\\
&\qquad + o\left(\sum_{i=m_h}^{n-2}|\mathbf{x}_i'(k)V_{i-1}\mathbf{x}_{i-1}(k)|\right)+O(1) \qquad \text{a.s.}
\end{aligned} \tag{3.32}
$$

By an argument similar to that used for showing Lemma 2.1 of Wei (1992), the Cauchy–Schwarz inequality and (3.13), we have

$$
\sum_{i=m_h}^{n-2}\mathbf{x}_i'(k)V_{i-1}\mathbf{x}_{i-1}(k)c_i^{-1} = \operatorname{tr}(\Gamma^{-1}(k)E_1(k))\log n + o(\log n) \qquad \text{a.s.},
$$

where $E_1(k) = E(\mathbf{x}_k(k)\mathbf{x}_{k+1}'(k))$, and

$$
\sum_{i=m_h}^{n-2}|\mathbf{x}_i'(k)V_{i-1}\mathbf{x}_{i-1}(k)| = O(\log n) \qquad \text{a.s.}
$$

These results, (3.32) and the fact that $\operatorname{tr}(\Gamma^{-1}(k)E_1(k)) = a_1(1,k)$ [note that $a_1(1,k) = a_1$ as $k \geq p_1$; see Section 1 for the definition of $a_j(h,k)$] together imply that

$$
F(k) = a_1 a_1(1,k)\sigma^2\log n + o(\log n) \qquad \text{a.s.} \tag{3.33}
$$

In view of (3.27)–(3.31) and (3.33) we have

$$\begin{aligned} C(k) &= 2a_1 a_1(1,k)\sigma^2 \log n \\ &\quad + o\left(\sum_{i=m_h}^{n-2}\left\{\mathbf{x}_i'(k)V_{i-1}\left(\sum_{j=k}^{i-2}\mathbf{x}_j(k)\eta_{j,2}\right)\right\}^2\right) + o(\log n) \qquad \text{a.s.} \end{aligned} \tag{3.34}$$

Since

$$\begin{aligned} &\sum_{i=m_h}^{n-2}\left\{\mathbf{x}_i'(k)V_{i-1}\left(\sum_{j=k}^{i-2}\mathbf{x}_j(k)\eta_{j,2}\right)\right\}^2 \\ &\qquad = \sum_{i=m_h}^{n-2}\left\{\mathbf{x}_i'(k)V_{i-1}\left(\sum_{j=k}^{i-1}\mathbf{x}_j(k)\eta_{j,2} - \mathbf{x}_{i-1}(k)\eta_{i-1,2}\right)\right\}^2, \end{aligned} \tag{3.35}$$

by the Cauchy–Schwarz inequality and an argument similar to that used to show (3.29), the right-hand side of (3.35) equals

$$(1+o(1))T(k) + o(\log n) \qquad \text{a.s.} \tag{3.36}$$

This fact, (3.22), (3.23), (3.26) and (3.34) yield

$$(1+o(1))T(k) = \{(1+a_1^2)k + 2a_1 a_1(1,k)\}\sigma^2 \log n + o(\log n) \qquad \text{a.s.} \tag{3.37}$$

According to (2.3), (3.21) and (3.37), (3.20) is obtained if we can show that

$$\sum_{i=m_h}^{n-2}\left\{\mathbf{x}_i'(k)V_{i-2}\left(\sum_{j=k}^{i-2}\mathbf{x}_j(k)\eta_{j,2}\right)\right\}^2 = T(k) + o(\log n) \qquad \text{a.s.} \tag{3.38}$$

To show (3.38), first observe that

$$\begin{aligned} &\mathbf{x}_i'(k)V_{i-1}\left(\sum_{j=1}^{i-1}\mathbf{x}_j(k)\eta_{j,2}\right) \\ &\qquad = \mathbf{x}_i'(k)V_{i-2}\sum_{j=k}^{i-2}\mathbf{x}_j(k)\eta_{j,2} + \mathbf{x}_i'(k)V_{i-2}\mathbf{x}_{i-1}(k)\eta_{i-1,2} \\ &\qquad\quad - \frac{\mathbf{x}_i'(k)V_{i-2}\mathbf{x}_{i-1}(k)}{1+\mathbf{x}_{i-1}'(k)V_{i-2}\mathbf{x}_{i-1}(k)}\mathbf{x}_{i-1}'(k)V_{i-2}\sum_{j=k}^{i-2}\mathbf{x}_j(k)\eta_{j,2} \\ &\qquad\quad - \frac{\mathbf{x}_i'(k)V_{i-2}\mathbf{x}_{i-1}(k)}{1+\mathbf{x}_{i-1}'(k)V_{i-1}\mathbf{x}_{i-1}(k)}\mathbf{x}_{i-1}'(k)V_{i-2}\mathbf{x}_{i-1}(k)\eta_{i-1,2}. \end{aligned}$$

This fact, Theorem 4 of Lai and Wei (1983), and an argument similar to that used to show (3.36) yield

$$T(k) = (1+o(1))\sum_{i=m_h}^{n-2}\left(\mathbf{x}_i'(k)V_{i-2}\sum_{j=k}^{i-2}\mathbf{x}_j(k)\eta_{j,2}\right)^2 + o(\log n) \qquad \text{a.s.,}$$

as asserted. □

Remark 3. Interestingly, it can be seen from Corollary 2.4 and Theorems 2.1, 2.2, 3.1 and 3.2 that the constant associated with the $1/n$ term of $\mathrm{MSPE}\, P_{n,h}(k)$, $f_{1,h}(k)$, appears in the $\log n$ term of $\mathrm{APE}\, P_{n,h}(k)$ and that associated with the $1/n$ term of $\mathrm{MSPE}\, D_{n,h}(k)$, $f_{2,h}(k)$, appears in the $\log n$ term of $\mathrm{APE}\, D_{n,h}(k)$. When $p_1$ and $p_h$ are known, these special features allow determination of the sign of $f_{1,h}(p_1) - f_{2,h}(p_h)$ by comparing the values of $\mathrm{APE}\, P_{n,h}(p_1)$ and $\mathrm{APE}\, D_{n,h}(p_h)$. This is because, according to (3.7) and (3.20), if $f_{1,h}(p_1) > f_{2,h}(p_h)$, then

$$P(\mathrm{APE}\, P_{n,h}(p_1) > \mathrm{APE}\, D_{n,h}(p_h) \text{ eventually}) = 1 \tag{3.39}$$

and if $f_{1,h}(p_1) < f_{2,h}(p_h)$, then

$$P(\mathrm{APE}\, P_{n,h}(p_1) < \mathrm{APE}\, D_{n,h}(p_h) \text{ eventually}) = 1. \tag{3.40}$$

Equalities (3.39) and (3.40) show that if $f_{1,h}(p_1) \neq f_{2,h}(p_h)$, then with probability 1 the sign of $\mathrm{APE}\, P_{n,h}(p_1) - \mathrm{APE}\, D_{n,h}(p_h)$ ultimately equals the sign of $f_{1,h}(p_1) - f_{2,h}(p_h)$.

Theorem 3.3 below deals with the asymptotic performances of $\mathrm{APE}\, P_{n,h}(k)$ and $\mathrm{APE}\, D_{n,h}(k)$ in underspecified cases.

Theorem 3.3. *Let the assumptions of Theorem 3.1 hold. Then for $1 \le k < p_1$ and $h \ge 1$,*

$$\begin{aligned} &\frac{1}{n}\left(\mathrm{APE}\, P_{n,h}(k) - \sum_{i=m_h}^{n-h} \eta_{i,h}^2\right) \\ &\quad = (\mathbf{a}_D(h,p_1) - \mathbf{a}_D(h,k))'\Gamma(p_1)(\mathbf{a}_D(h,p_1) - \mathbf{a}_D(h,k)) \\ &\qquad + (\mathbf{a}(h,k) - \mathbf{a}_D(h,k))'\Gamma(k)(\mathbf{a}(h,k) - \mathbf{a}_D(h,k)) \\ &\qquad + o(1) \qquad a.s., \end{aligned} \tag{3.41}$$

*where $\mathbf{a}(h,k) = A^{h-1}(k)\mathbf{a}_D(1,k)$ with $A(k)$ defined after (2.2) and $\mathbf{a}_D(h,k)$ in the first term of the right-hand side viewed as a $p_1$-dimensional vector with undefined entries set to zero, and for $1 \le k < p_h$ and $h \ge 1$,*

$$\begin{aligned} &\frac{1}{n}\left(\mathrm{APE}\, D_{n,h}(k) - \sum_{i=m_h}^{n-h} \eta_{i,h}^2\right) \\ &\quad = (\mathbf{a}_D(h,p_h) - \mathbf{a}_D(h,k))'\Gamma(p_h)(\mathbf{a}_D(h,p_h) - \mathbf{a}_D(h,k)) \\ &\qquad + o(1) \qquad a.s., \end{aligned} \tag{3.42}$$

*where $\mathbf{a}_D(h,k)$ in the right-hand side is viewed as a $p_h$-dimensional vector with undefined entries set to zero.*

PROOF. Following Hemerly and Davis (1989) [which deals with APE $P_{n,h}(k)$ with $h=1$], we have

$$\begin{aligned}\text{APE}\, P_{n,h}(k) &= \sum_{i=m_h}^{n-h}\{\eta_{i,h} + \mathbf{x}_i'(p_1)(\mathbf{a}_D(h,p_1) - \hat{\mathbf{a}}_i(h,k))\}^2 \\ &= \sum_{i=m_h}^{n-h}\eta_{i,h}^2 + (1+o(1))\sum_{i=m_h}^{n-h}\{\mathbf{x}_i'(p_1)(\mathbf{a}_D(h,p_1) - \hat{\mathbf{a}}_i(h,k))\}^2 \\ &\quad + O(1) \qquad \text{a.s.,}\end{aligned} \tag{3.43}$$

where $\hat{\mathbf{a}}_i(h,k)$ is now viewed as a $p_1$-dimensional vector with undefined entries set to zero and the second equality is ensured by Chow (1965). Since (3.11) ensures that $\lim_{n\to\infty}\hat{\mathbf{a}}_n(h,k) = \mathbf{a}(h,k)$ a.s., we can rewrite (3.43) as

$$\begin{aligned}\text{APE}\, P_{n,h}(k) &= (1+o(1))(\mathbf{a}_D(h,p_1) - \mathbf{a}(h,k))' \\ &\quad \times \sum_{i=m_h}^{n-h}\mathbf{x}_i(p_1)\mathbf{x}_i'(p_1)(\mathbf{a}_D(h,p_1) - \mathbf{a}(h,k)) \\ &\quad + \sum_{i=m_h}^{n-h}\eta_{i,h}^2 + o\left(\sum_{i=m_h}^{n-h}\mathbf{x}_i'(p_1)\mathbf{x}_i(p_1)\right) + O(1) \qquad \text{a.s.}\end{aligned}$$

Consequently, (3.41) follows from (3.11) and the fact that

$$(\mathbf{a}_D(h,p_1) - \mathbf{a}_D(h,k))'\Gamma(p_1)(\mathbf{a}_D(h,k) - \mathbf{a}(h,k)) = 0,$$

where $\mathbf{a}_D(h,k)$ and $\mathbf{a}(h,k)$ are viewed as $p_1$-dimensional vectors with undefined entries set to zero.

Since the proof for (3.42) is similar to that for (3.41), to save space we omit the details. □

Armed with the previous results, we are now in a position to show the asymptotic efficiency of $(\hat{k}_n, \hat{j}_n)$.

THEOREM 3.4. *Let the assumptions of Theorem* 3.1 *hold. Then, for* $K \geq p_1$ $(\hat{k}_n \hat{j}_n)$ *is asymptotically efficient in the sense of* (3.3).

PROOF. First note that for $k \geq p_1$, $f_{2,1}(k) = k$. Hence Theorem 3.2 yields that for $k > p_1$, $P(\text{APE}\, D_{n,1}(p_1) < \text{APE}\, D_{n,1}(k) \text{ eventually}) = 1$. Since the first term on the right-hand side of (3.42) is positive, by Theorems 3.2 and 3.3 we have for $k < p_1$, $P(\text{APE}\, D_{n,1}(p_1) < \text{APE}\, D_{n,1}(k) \text{ eventually}) = 1$. As a result, $\hat{k}_{D,n}^{(1)} = p_1 + o(1)$ a.s. This fact and Theorems 3.1–3.3 further ensure that

$$P((\hat{k}_n, \hat{j}_n) \in C_{h,K} \text{ eventually}) = 1,$$

as asserted. □

REMARK 4. In this remark, we consider the problem of choosing $p_h$, $h \geq 1$, under model (1.1). For $h = 1$ we have shown in the proof of Theorem 3.4 that

$$\hat{k}_{D,n}^{(h)} = p_h + o(1) \qquad \text{a.s.} \tag{3.44}$$

This motivated us to ask whether (3.44) still holds with $h \geq 2$. To investigate this question first assume $p_h = p_1$ (or, equivalently, $b_{h-1} \neq 0$). By (ii) of Theorem 2.3 and Theorems 3.2 and 3.3, this assumption guarantees that (3.44) holds with $h \geq 2$. [In fact, by (i) of Theorem 2.3 and Theorems 3.1 and 3.3, this assumption also ensures that for $h \geq 2$,

$$\lim_{n\to\infty} \hat{k}_{P,n}^{(h)} = p_1 = p_h \qquad \text{a.s.},$$

where $\hat{k}_{P,n}^{(h)} = \arg\min_{1\leq k\leq K} \mathrm{APE}\, P_{n,h}(k)$.] However, when $h$ is large and $p_h \leq k < p_1$ it is very difficult to verify $f_{2,h}(k) < f_{2,h}(k+1)$, which is an essential property for (3.44) with $h \geq 2$ to be true. [Note that (2.10) only ensures that $f_{2,h}(k) < f_{2,h}(k+1)$ holds with $k \geq p_1$.] Consequently, with arguments used in the present article, (3.44) cannot be guaranteed without extra constraints on the parameter space.

To establish a strongly consistent estimator of $p_h$ without constraints on the parameter space, we consider the multistep generalization of the Bayesian information criterion (BIC),

$$\mathrm{BIC}_{n,h}(k) = \log \hat{\sigma}^2_{D_h,n}(k) + \frac{kc_n}{n},$$

where $h \geq 1, c_n \to \infty, c_n = o(n)$, $\liminf_{n\to\infty} c_n/(\log n) > 0$ and $\hat{\sigma}^2_{D_h,n}(k) = (1/n)\sum_{i=k}^{n-h}(x_{i+h} - \mathbf{x}_i'(k)\check{\mathbf{a}}_n(h,k))^2$. When the assumptions of Theorem 3.2 hold, then arguments similar to those used to show Theorem 3.2 of the present study and Theorem 3.6 of Wei (1992) yield that

$$\hat{k}_{B,n}^{(h)} = p_h + o(1) \qquad \text{a.s.},$$

where $\hat{k}_{B,n}^{(h)} = \arg\min_{1\leq k\leq K} \mathrm{BIC}_{n,h}(k)$. Therefore, the difficulty encountered with $\hat{k}_{D,n}^{(h)}$ does not exist for $\hat{k}_{B,n}^{(h)}$.

**4. An extension to subset autoregressions.** When some $a_i$'s with $1 \leq i \leq p_1 - 1$ in model (1.1) or some $a_i(h, p_h)$'s with $1 \leq i \leq p_h - 1$ in model (1.7) are zero, a multistep predictor, which is obtained without estimating these zero coefficients, can be more efficient than the best predictor among families I and II. This motivated us to consider the selection of subset autoregressive models. Several different algorithms are available for choosing the one-step prediction model under this more general setting [e.g., McClave (1975) and Haggan and Oyetunji (1984)]. While these algorithms have their own advantages, no algorithm has been shown to possess optimal properties from the (multistep) MSPE point of view. An algorithm which is modified from $(\hat{k}_n, \hat{j}_n)$ is therefore proposed in this section as a remedy.

To begin with, let $\theta_i = 1$ if $x_{t+1-i}$ is included as a regressor variable for predicting $x_{t+h}$ and let $\theta_i = 0$ if $x_{t+1-i}$ is not included. Then the family of all (nontrivial) subset autoregressions can be expressed as

$$\Theta = \{\theta = (\theta_1, \ldots, \theta_K) : \theta_i = 0 \text{ or } 1 \text{ for } 1 \le i \le K, \text{ and } \theta_i = 1 \text{ for at least one } i\},$$

where $K$ is as defined in Section 1. When model $\theta \in \Theta$ is adopted, the corresponding plug-in and direct predictors of $x_{n+h}$ are denoted by $\hat{x}_{n+h}(\theta)$ [or $(\theta, 1)$] and $\check{x}_{n+h}(\theta)$ [or $(\theta, 2)$], respectively, and the multistep MSPEs of $\hat{x}_{n+h}(\theta)$ and $\check{x}_{n+h}(\theta)$ are denoted by $\mathrm{MSPE}\, P_{n,h}(\theta)$ and $\mathrm{MSPE}\, D_{n,h}(\theta)$, respectively. In addition, we also use $\mathrm{APE}\, P_{n,h}(\theta)$ and $\mathrm{APE}\, D_{n,h}(\theta)$, respectively, to denote the multistep APEs based on sequential plug-in and direct predictors when $\theta \in \Theta$ is used. Let $\theta^{(1)} = (\theta_1^{(1)}, \ldots, \theta_K^{(1)})$ and $\theta^{(2)} = (\theta_1^{(2)}, \ldots, \theta_K^{(2)})$ be members of $\Theta$. Then we say $\theta^{(1)} \le \theta^{(2)}$ if $\theta_i^{(1)} \le \theta_i^{(2)}$ for all $1 \le i \le K$ and $\theta^{(1)} \not\le \theta^{(2)}$ if $\theta_i^{(1)} > \theta_i^{(2)}$ for at least one $i$. Now the modified model selection procedure $(\hat{\theta}_n\ \hat{j}_n)$ with $\hat{\theta}_n \in \Theta$ and $1 \le \hat{j}_n \le 2$, is given as follows.

STEP 1. Define $\hat{\theta}_{D,n}^{(1)} = \arg\min_{\theta \in \Theta} \mathrm{APE}\, D_{n,1}(\theta)$.

STEP 2. Define

$$\hat{\theta}_{D,n}^{(h)} = \arg\min_{\theta \in \Theta} \mathrm{APE}\, D_{n,h}(\theta)$$

and define

$$\hat{\theta}_n^{(1,h)} = \arg\min_{\theta \in \Theta_1} \mathrm{APE}\, P_{n,h}(\theta),$$

where $\Theta_1 = \{\theta : \theta \in \Theta \text{ and } \hat{\theta}_{D,n}^{(1)} \le \theta\}$.

STEP 3. If $\mathrm{APE}\, D_{n,h}(\hat{\theta}_{D,n}^{(h)}) > \mathrm{APE}\, P_{n,h}(\hat{\theta}_n^{(1,h)})$, then $(\hat{\theta}_n, \hat{j}_n) = (\hat{\theta}_n^{(1,h)}, 1)$; otherwise $(\hat{\theta}_n, \hat{j}_n) = (\hat{k}_{D,n}^{(h)}, 2)$.

To show the validity of $(\hat{\theta}_n, \hat{j}_n)$, let us recall models (1.1) and (1.7) again, and define $\theta^* = (\theta_1^*, \ldots, \theta_K^*)$ and $\theta^{**} = (\theta_1^{**}, \ldots, \theta_K^{**})$, where $\theta_i^* = 1$ if $a_i \ne 0$ and $\theta_i^* = 0$ if $a_i = 0$ or $i > p_1$, and $\theta_i^{**} = 1$ if $a_i(h, p_h) \ne 0$ and $\theta_i^{**} = 0$ if $a_i(h, p_h) = 0$ or $i > p_h$. Therefore, $\theta^*$ and $\theta^{**}$, respectively, are the most parsimonious correct models for the plug-in and direct predictors. Following (3.1) and (3.2), the loss functions of $\hat{x}_{n+h}(\theta)$ and $\check{x}_{n+h}(\theta)$ are defined as

$$E_{1,h}(\theta) = \begin{cases} \lim_{n\to\infty} n(\text{MSPE}\, P_{n,h}(\theta) - \sigma_h^2), & \text{if } \theta^* \le \theta, \\ \infty, & \text{if } \theta^* \not\le \theta, \end{cases} \tag{4.1}$$

and

$$E_{2,h}(k) = \begin{cases} \lim_{n\to\infty} n(\text{MSPE}\, D_{n,h}(\theta) - \sigma_h^2), & \text{if } \theta^{**} \le \theta, \\ \infty, & \text{if } \theta^{**} \not\le \theta, \end{cases} \tag{4.2}$$

respectively, where the existence of the above limits is guaranteed by arguments similar to those used to obtain Theorems 2.1 and 2.2. [Note that we also obtain expressions for the above limits like those on the right-hand sides of (2.2) and (2.3). However, these expressions are not presented here, since they are not needed in the following analysis.] A model selection criterion $(\tilde{\theta}_n, \tilde{j}_n)$ with $\tilde{\theta}_n \in \Theta$ and $1 \le \tilde{j}_n \le 2$ is said to be asymptotically efficient if

$$P((\tilde{\theta}_n, \tilde{j}_n) \in B_{h,K} \text{ eventually}) = 1, \tag{4.3}$$

where

$$B_{h,K} = \left\{(\theta, j) : \theta \in \Theta,\ 1 \le j \le 2 \text{ and } E_{j,h}(\theta) = \min_{\theta_0 \in \Theta, 1 \le j_0 \le 2} E_{j_0,h}(\theta_0)\right\}.$$

The main result of this section is stated as follows.

THEOREM 4.1. *Let the assumptions of Theorem* 3.1 *hold. Then* $(\hat{\theta}_n, \hat{j}_n)$ *is asymptotically efficient in the sense of* (4.3).

Theorem 4.1 can be shown by arguments similar to those used to show Theorems 3.1–3.4. To save space, the details are omitted. Theorems 3.4 and 4.1 yield that for sufficiently large $n$, the predictor selected by $(\hat{\theta}_n, \hat{j}_n)$ is at least as efficient as the one selected by $(\hat{k}_n, \hat{j}_n)$. Before leaving this section, we note that the main disadvantage of $(\hat{\theta}_n, \hat{j}_n)$ is its time-consuming nature, since it needs to compute the multistep APEs for all possible subset autoregressive models and for two different prediction methods. However, with the availability of fast computers and efficient recursive formulas the computer time needed to complete this task is not expensive, provided $K$ is not too large.

**5. Concluding remarks.** One of the main purposes of this article was to find the optimal multistep predictor in finite-order AR models from the honest MSPE point of view. Since both the plug-in and the direct predictors are considered, it is not possible to achieve this goal by identifying the order of the smallest correct model, as discussed in Section 2. To resolve this problem, a new predictor selection procedure, $(\hat{k}_n, \hat{j}_n)$ is proposed. We show that for sufficiently large $n$, $(\hat{k}_n, \hat{j}_n)$ can achieve the above goal by choosing the best combination of the prediction order and the prediction method. In Section 4 this procedure is extended to the situation where all possible subset autoregressions are included as candidate models. On the other hand, the parameter set where (1.8) occurs has Lebesgue measure zero. So one may argue that this is unlikely to occur in practice and, hence, the necessity to construct $(\hat{k}_n, \hat{j}_n)$ may be questioned. In contrast to this criticism, it is worth noting that $(\hat{k}_n, \hat{j}_n)$ asymptotically dominates traditional multistep prediction procedures, which select the one-step prediction order by certain consistent order selection criteria and then forecast $x_{n+h}$ through the plug-in (or direct) method. More precisely, the predictor selected by $(\hat{k}_n, \hat{j}_n)$ has at least the same asymptotic efficiency as those predictors selected by the traditional procedures for all points of $\Lambda$ and is asymptotically more efficient than the latter for some nonempty subset of $\Lambda$ [since the set where (1.8) occurs is nonempty for $h \geq 2$]. Moreover, some other advantages of $(\hat{k}_n, \hat{j}_n)$, besides offering a treatment of the case where (1.8) occurs, are also emphasized at the end of Section 1.

The validity of $(\hat{k}_n, \hat{j}_n)$ is justified in the stationary case. It is also believed that the predictor chosen by this procedure may also perform well in unstable cases. However, since the proofs of Theorems 3.1 and 3.2 (especially Theorem 3.1) rely highly on stationary assumptions, their extensions to unstable cases are not straightforward. Further work is needed to overcome these technical difficulties.

This article assumes that the order of the underlying AR model is finite. Hence, the frequently discussed AR($\infty$) model is excluded. When the data are known to be generated from an AR($\infty$) model, it is common to use an AR model of increasing (with $n$) order to predict future observations; see, for example, Shibata (1980), Gerencsér (1992), Bhansali (1996) and Ing and Wei (2003, 2004). In this situation, Ing and Wei (2004) showed that AIC is asymptotically efficient for the honest one-step prediction. On the other hand, Ing and Yu (2002) showed that the one-step APE is not asymptotically efficient in this situation. To rectify the difficulty of using APE in AR($\infty$) models, Ing and Yu (2002) proposed a modification of APE, $\mathrm{APE}_\delta$. Instead of accumulating squares of sequential prediction errors from stage $m_1$ [see (1.9)], $\mathrm{APE}_\delta$ is obtained by accumulating squares of sequential prediction errors from stage $n\delta$, where $0 < \delta < 1$ may depend on $n$. Under certain regularity conditions, they showed that $\mathrm{APE}_\delta$ is asymptotically efficient in AR($\infty$) models. Motivated by this result, it is expected that an efficient multistep predictor selection criterion can be established in an AR($\infty$) model after asymptotic behavior of $\mathrm{APE}\, P_{n,h}(k)$ and of $\mathrm{APE}\, D_{n,h}(k)$, with $h \geq 2$ and $m_h$ replaced by $n\delta$, $0 < \delta < 1$, is clarified under this model. As a final remark, we note that when it is a priori unknown whether the order of the underlying AR model is finite or infinite, the choice between the original APE and its modification (by Ing and Yu) becomes a challenging problem even for one-step predictions. Can a modification of $(\hat{k}_n, \hat{j}_n)$ be obtained for the optimal multistep prediction without order assumptions? This is the subject of ongoing research.

## APPENDIX

Proof of (3.16). By (3.11), Theorem 3 of Lai and Wei (1983) and Chow (1965), (3.16) is guaranteed by showing that

$$\sum_{i=m_h}^{n-h} \|\mathbf{x}_i(k)\|^2 \left\| \frac{1}{i-k} \sum_{j=k}^{i-1} \mathbf{x}_j(k)\varepsilon_{j+1} \right\|^2 = O(\log n) \qquad \text{a.s.} \tag{A.1}$$

To obtain (A.1), first observe that the term on the left-hand side of (A.1) can be expressed as

$$\begin{aligned} &\sum_{i=m_h}^{n-h} \left\{ \left( \sum_{l=0}^{k-1} x_{i-l}^2 \right) \left( \frac{1}{(i-k)^2} \sum_{j_1=k}^{i-1} \sum_{j_2=k}^{i-1} \left( \sum_{c=0}^{k-1} x_{j_1-c} x_{j_2-c} \right) \varepsilon_{j_1+1}\varepsilon_{j_2+1} \right) \right\} \\ &\quad = \sum_{l=0}^{k-1} \sum_{c=0}^{k-1} \left\{ \sum_{i=m_h}^{n-h} \left( \frac{1}{(i-k)^2} \sum_{j_1=k}^{i-1} \sum_{j_2=k}^{i-1} x_{j_1-c} x_{j_2-c} \varepsilon_{j_1+1}\varepsilon_{j_2+1} \right) x_{i-l}^2 \right\}. \end{aligned} \tag{A.2}$$

In view of (A.2), if we can show that

$$\sum_{i=m_h}^{n-h} \left( \frac{1}{(i-k)^2} \sum_{j_1=k}^{i-1} \sum_{j_2=k}^{i-1} x_{j_1-c} x_{j_2-c} \varepsilon_{j_1+1}\varepsilon_{j_2+1} \right) x_{i-l}^2 = O(\log n) \qquad \text{a.s.} \tag{A.3}$$

for each $0 \le l \le k-1$ and $0 \le c \le k-1$, then (A.1) follows. In what follows we prove this property only for the case of $c = l = 0$, because the results for other $c$'s and $l$'s can be obtained similarly.

Note that

$$\begin{aligned} &\sum_{i=m_h}^{n-h} \left( \frac{1}{(i-k)^2} \sum_{j_1=k}^{i-1} \sum_{j_2=k}^{i-1} x_{j_1} x_{j_2} \varepsilon_{j_1+1}\varepsilon_{j_2+1} \right) x_i^2 \\ &\quad \le C^* \sum_{i=k+1}^{n-h} \left( \frac{1}{i^2} \sum_{j_1=k}^{i-1} \sum_{j_2=k}^{i-1} x_{j_1} x_{j_2} \varepsilon_{j_1+1}\varepsilon_{j_2+1} \right) x_i^2 \\ &\quad = C^* \sum_{j_1=k}^{n-h-1} \sum_{j_2=k}^{n-h-1} x_{j_1} x_{j_2} \varepsilon_{j_1+1}\varepsilon_{j_2+1} \left( \sum_{i=r}^{n-h} \frac{x_i^2}{i^2} \right), \end{aligned} \tag{A.4}$$

where $C^*$ is some positive number and $r = \max\{j_1+1, j_2+1\}$. Observe that

$$\begin{aligned}\sum_{i=r}^{n-h}\frac{x_i^2}{i^2} &= \sum_{i=r}^{n-h}\left(\frac{s_i^2 - i\gamma_0}{i^2}\right) - \left(\frac{s_{i-1}^2-(i-1)\gamma_0}{(i-1)^2}\right)\\ &\quad + \sum_{i=r}^{n-h}\frac{(s_{i-1}^2-(i-1)\gamma_0)(2i-1)}{(i-1)^2 i^2} + \gamma_0\sum_{i=r}^{n-h}\frac{1}{i^2}\\ &= A_n + B_{n,r} + C_{n,r} + D_{n,r},\end{aligned}$$

where $s_i^2 = \sum_{j=1}^{i} x_j^2$, $\gamma_0 = E(x_1^2)$,

$$\begin{aligned}A_n &= \frac{s_{n-h}^2 - (n-h)\gamma_0}{(n-h)^2},\\ B_{n,r} &= -\frac{s_{r-1}^2-(r-1)\gamma_0}{(r-1)^2},\\ C_{n,r} &= \sum_{i=r}^{n-h}\frac{(s_{i-1}^2-(i-1)\gamma_0)(2i-1)}{(i-1)^2 i^2}\end{aligned}$$

and

$$D_{n,r} = \gamma_0\sum_{i=r}^{n-h} i^{-2}.$$

This and (A.4) yield

$$\begin{aligned}&\sum_{i=k+1}^{n-h}\left(\frac{1}{i^2}\sum_{j_1=k}^{i-1}\sum_{j_2=k}^{i-1}x_{j_1}x_{j_2}\varepsilon_{j_1+1}\varepsilon_{j_2+1}\right)x_i^2\\ &\quad= \sum_{j_1=k}^{n-h-1}\sum_{j_2=k}^{n-h-1}x_{j_1}x_{j_2}\varepsilon_{j_1+1}\varepsilon_{j_2+1}(A_n+B_{n,r}+C_{n,r}+D_{n,r}).\end{aligned}\tag{A.5}$$

Since

$$\sum_{j_1=k}^{n-h-1}\sum_{j_2=k}^{n-h-1}x_{j_1}x_{j_2}\varepsilon_{j_1+1}\varepsilon_{j_2+1}A_n = o(1)\frac{1}{n}\left(\sum_{j=k}^{n-h-1}x_j\varepsilon_{j+1}\right)^2 \qquad \text{a.s.},$$

by Wei [(1987), equation (2.30)] and (3.11),

$$\sum_{j_1=k}^{n-h-1}\sum_{j_2=k}^{n-h-1}x_{j_1}x_{j_2}\varepsilon_{j_1+1}\varepsilon_{j_2+1}A_n = o(\log n) \qquad \text{a.s.}\tag{A.6}$$

By (3.11), an analogy with Lemma 2.1 of Wei (1992) and Chow (1965),

$$
\begin{aligned}
&\sum_{j_1=k}^{n-h-1}\sum_{j_2=k}^{n-h-1} x_{j_1}x_{j_2}\varepsilon_{j_1+1}\varepsilon_{j_2+1}B_{n,r}\\
&\quad = -\biggl\{\sum_{j=k}^{n-h-1}\frac{x_j^2\varepsilon_j^2(s_{j-1}^2-(j-1)\gamma_0)}{(j-1)^2}\\
&\qquad\quad + 2\sum_{j_2=k+1}^{n-h-1}\Biggl(\sum_{j_1=k}^{j_2-1}x_{j_1}\varepsilon_{j_1+1}\Biggr)\frac{s_{j_2-1}^2-(j_2-1)\gamma_0}{(j_2-1)^2}x_{j_2}\varepsilon_{j_2+1}\biggr\}\\
&\quad = o(\log n)+o\Biggl(\sum_{j_2=k+1}^{n-h-1}\Biggl(\frac{1}{j_2}\sum_{j_1=k}^{j_2-1}x_{j_1}\varepsilon_{j_1+1}\Biggr)^2x_{j_2}^2\Biggr)\qquad \text{a.s.}
\end{aligned}
\tag{A.7}
$$

Exchanging the order of summation, we have

$$
\begin{aligned}
&\sum_{j_1=k}^{n-h-1}\sum_{j_2=k}^{n-h-1} x_{j_1}x_{j_2}\varepsilon_{j_1+1}\varepsilon_{j_2+1}C_{n,r}\\
&\quad = \sum_{i=k+1}^{n-h}\Biggl(\sum_{j=k}^{i-1}x_j\varepsilon_{j+1}\Biggr)^2\frac{\{s_{i-1}^2-(i-1)\gamma_0\}(2i-1)}{i^2(i-1)^2}\\
&\quad = o\Biggl(\sum_{i=k+1}^{n-h}\Biggl(\sum_{j=k}^{i-1}x_j\varepsilon_{j+1}\Biggr)^2\frac{1}{i^2}\Biggr)\qquad \text{a.s.,}
\end{aligned}
\tag{A.8}
$$

where the second equality is ensured by (3.11). Observe that

$$
\begin{aligned}
&\sum_{i=k+1}^{n-h}\Biggl(\sum_{j=k}^{i-1}x_j\varepsilon_{j+1}\Biggr)^2\frac{1}{i^2}\\
&\quad = \sum_{j_1=k}^{n-h-1}\sum_{j_2=k}^{n-h-1}x_{j_1}x_{j_2}\varepsilon_{j_1+1}\varepsilon_{j_2+1}\sum_{i=r}^{n-h}\frac{1}{i^2}\\
&\quad = \sum_{j=k}^{n-h-1}x_j^2\varepsilon_{j+1}^2\sum_{i=j+1}^{n-h}\frac{1}{i^2}\\
&\qquad + 2\sum_{j_2=k+1}^{n-h-1}\Biggl(\sum_{j_1=k}^{j_2-1}x_{j_1}\varepsilon_{j_1+1}\Biggr)x_{j_2}\Biggl(\sum_{i=j_2+1}^{n-h}\frac{1}{i^2}\Biggr)\varepsilon_{j_2+1}\\
&\quad = O(\log n)+o\Biggl(\sum_{j_2=k+1}^{n-h-1}\Biggl(\frac{1}{j_2}\sum_{j_1=k}^{j_2-1}x_{j_1}\varepsilon_{j_1+1}\Biggr)^2x_{j_2}^2\Biggr)\qquad \text{a.s.,}
\end{aligned}
\tag{A.9}
$$

where the last equality follows from an argument similar to that used for showing (A.7). As a result, (A.8) and (A.9) yield

$$\begin{aligned}\text{(A.10)}\quad &\sum_{j_1=k}^{n-h-1}\sum_{j_2=k}^{n-h-1} x_{j_1}x_{j_2}\varepsilon_{j_1+1}\varepsilon_{j_2+1}C_{n,r}\\ &\quad = o(\log n) + o\left(\sum_{j_2=k+1}^{n-h-1}\left(\frac{1}{j_2}\sum_{j_1=k}^{j_2-1} x_{j_1}\varepsilon_{j_1+1}\right)^2 x_{j_2}^2\right) \qquad \text{a.s.}\end{aligned}$$

Reasoning as for (A.9),

$$\begin{aligned}\text{(A.11)}\quad &\sum_{j_1=k}^{n-h-1}\sum_{j_2=k}^{n-h-1} x_{j_1}x_{j_2}\varepsilon_{j_1+1}\varepsilon_{j_2+1}D_{n,r}\\ &\quad = o\left(\sum_{j_2=k+1}^{n-h-1}\left(\frac{1}{j_2}\sum_{j_1=k}^{j_2-1} x_{j_1}\varepsilon_{j_1+1}\right)^2 x_{j_2}^2\right) + O(\log n) \qquad \text{a.s.}\end{aligned}$$

Consequently, (A.3) [and hence (A.1)] follows from (A.4)–(A.7), (A.10) and (A.11). □

**Acknowledgments.** I am deeply grateful to Dr. David Findley and Professor Ching Zong Wei for valuable discussions and suggestions on a previous version of this paper. I also thank an Editor, an Associate Editor and three anonymous referees for helpful comments that substantially improved the presentation.

## REFERENCES

AKAIKE, H. (1969). Fitting autoregressive models for prediction. *Ann. Inst. Statist. Math.* **21** 243–247. MR246476

AKAIKE, H. (1974). A new look at the statistical model identification. *IEEE Trans. Automat. Control* **19** 716–723. MR423716

BHANSALI, R. J. (1996). Asymptotically efficient autoregressive model selection for multistep prediction. *Ann. Inst. Statist. Math.* **48** 577–602. MR1424784

BHANSALI, R. J. (1997). Direct autoregressive predictors for multistep prediction: Order selection and performance relative to the plug in predictors. *Statist. Sinica* **7** 425–449. MR1466690

CHOW, Y. S. (1965). Local convergence of martingales and the law of large numbers. *Ann. Math. Statist.* **36** 552–558. MR182040

FINDLEY, D. F. (1984). On some ambiguities associated with the fitting of ARMA models to time series. *J. Time Ser. Anal.* **5** 213–225. MR782076

FINDLEY, D. F., PÖTSCHER, B. M. and WEI, C. Z. (2001). Uniform convergence of sample second moments of families of time series arrays. *Ann. Statist.* **29** 815–838. MR1865342

FINDLEY, D. F., PÖTSCHER, B. M. and WEI, C. Z. (2003). Modeling of time series arrays by multistep prediction or likelihood methods. *J. Econometrics* **118** 151–187. MR2030971

GERENCSÉR, L. (1992). AR($\infty$) estimation and nonparametric stochastic complexity. *IEEE Trans. Inform. Theory* **38** 1768–1778. MR1187818

Haggan, V. and Oyetunji, O. B. (1984). On the selection of subset autoregressive time series models. *J. Time Ser. Anal.* **5** 103–113.

Haywood, J. and Tunnicliffe-Wilson, G. (1997). Fitting time series models by minimising multistep-ahead errors: A frequency domain approach. *J. Roy. Statist. Soc. Ser. B* **59** 237–254. MR1436566

Hemerly, E. M. and Davis, M. H. A. (1989). Strong consistency of the PLS criterion for order determination of autoregressive processes. *Ann. Statist.* **17** 941–946. MR994279

Hurvich, C. M. and Tsai, C.-L. (1997). Selection of a multistep linear predictor for short time series. *Statist. Sinica* **7** 395–406. MR1466688

Ing, C.-K. (2003). Multistep prediction in autoregressive processes. *Econometric Theory* **19** 254–279. MR1966030

Ing, C.-K. and Wei, C. Z. (2003). On same-realization prediction in an infinite-order autoregressive process. *J. Multivariate Anal.* **85** 130–155. MR1978181

Ing, C.-K. and Wei, C. Z. (2004). Order selection for same-realization predictions in autoregressive processes. Submitted for publication. MR1978181

Ing, C.-K. and Yu, S. H. (2002). Asymptotic equivalence between information- and prediction-based model selection criteria: A new look at AIC. Unpublished manuscript.

Lai, T. L. and Wei, C. Z. (1982). Least squares estimates in stochastic regression models with applications to identification and control of dynamic systems. *Ann. Statist.* **10** 154–166. MR642726

Lai, T. L. and Wei, C. Z. (1983). Asymptotic properties of general autoregressive models and strong consistency of least squares estimates of their parameters. *J. Multivariate Anal.* **13** 1–23. MR695924

Lai, T. L. and Wei, C. Z. (1985). Asymptotic properties of multivariate weighted sums with applications to stochastic regression in linear dynamic systems. In *Multivariate Analysis VI* (P. R. Krishnaiah, ed.) 375–393. North-Holland, Amsterdam. MR822308

McClave, J. (1975). Subset autoregression. *Technometrics* **17** 213–220. MR368359

Rissanen, J. (1986). Order estimation by accumulated prediction errors. *J. Appl. Probab.* **23A** 55–61. MR803162

Rissanen, J. (1987). Stochastic complexity (with discussion). *J. Roy Statist. Soc. Ser. B* **49** 223–265. MR928936

Rissanen, J. (1989). *Stochastic Complexity in Statistical Inquiry*. World Scientific, New York. MR1082556

Schwarz, G. (1978). Estimating the dimension of a model. *Ann. Statist.* **6** 461–464. MR468014

Shibata, R. (1980). Asymptotically efficient selection of the order of the model for estimating parameters of a linear process. *Ann. Statist.* **8** 147–164. MR557560

Tiao, G. C. and Xu, D. (1993). Robustness of maximum likelihood estimates for multistep predictions: The exponential smoothing case. *Biometrika* **80** 623–641. MR1248027

Wei, C. Z. (1987). Adaptive prediction by least squares predictors in stochastic regression models with application to time series. *Ann. Statist.* **15** 1667–1682. MR913581

Wei, C. Z. (1992). On predictive least squares principles. *Ann. Statist.* **20** 1–42. MR1150333

Institute of Statistical Science
Academia Sinica
128 Academia road Section 2
Taipei 115
Taiwan
e-mail: cking@stat.sinica.edu.tw